\documentclass[twoside,a4paper,12pt]{article}
\usepackage{amsmath,amsthm,amsfonts,amssymb}
\usepackage{mathrsfs}
\usepackage{color,graphicx}

\makeatletter
\def\@seccntformat#1{\csname the#1\endcsname.\quad}
\renewcommand\section{\@startsection {section}{1}{\z@}%
                                   {-3.5ex \@plus -1ex \@minus -.2ex}%
                                   {2.3ex \@plus.2ex}%
                                   {\normalfont\bf\center }}
\renewcommand\subsection{\@startsection {subsection}{1}{\z@}%
                                   {-3.5ex \@plus -1ex \@minus -.2ex}%
                                   {2.3ex \@plus.2ex}%
                                   {\normalfont\bf}}
\makeatother

\date{\today}

\newcommand{\pa}{\partial}

\newtheoremstyle{new-thm}
 {3pt}
 {3pt}
 {\it}
 {0pt} 
 {\bf}
 {.}
 {.5em}
 {}
\newtheoremstyle{new-def}
 {3pt}
 {3pt}
 {\rm}
 {0pt} 
 {\bf}
 {.}
 {.5em}
 {}

\theoremstyle{new-thm}
  \newtheorem{thm}{Theorem}
  
  \newtheorem{lemma}[thm]{Lemma}

\theoremstyle{new-def}

  \newtheorem{rem}[thm]{Remark}

\begin{document}

\vspace*{2cm}
\begin{center} {\Large\bf On the zeros of the Macdonald functions}
\end{center} 

\bigskip

\begin{center} Yuji Hamana, Hiroyuki Matsumoto and 
Tomoyuki Shirai\end{center} 

\bigskip

\begin{quote} {\bf Abstract.} 
We are concerned with the zeros of the Macdonald functions 
or the modified Bessel functions of the second kind with real index. 
By using the explicit expressions for the algebraic equations 
satisfied by the zeros, 
we describe the behavior of the zeros when the index moves.   
Results by numerical computations are also presented.   

{\it keywords}: Zeros, Bessel functions, Asymptotic expansion
\end{quote}

In this article 
we are concerned with the zeros of the Macdonald function, 
or the modified Bessel function of the second kind with real index, 
which we denote by $K_\nu$ in the usual notation.   
By analytic continuation we consider $K_\nu(z)$ 
as a function in $z\in{\mathbf C}\setminus(-\infty,0]$.  
It is well known that $K_\nu(z)$ is an entire function in $\nu$.   
The zeros of the Bessel functions $J_\nu, Y_\nu$ and 
of the other modified function $I_\nu$ are well studied 
and we can also carry out the numerical computations 
for them in several ways.   
However, only a few things are known about the zeros of $K_\nu$.  
See \cite{KS,W,ZT}.   

Since $K_{-\nu}=K_\nu$ and 
there are no zeros when $0\leqq \nu<\frac{3}{2}$ (see \cite{W}), 
we throughout assume $\nu\geqq\frac{3}{2}$.   
When $\nu=2n+\frac{3}{2}, n=0,1,...,$ 
$K_{2n+\frac{3}{2}}$ is of the form 
$\sqrt{\pi/(2z)}z^{-(2n+1)}\varphi_n(z)$, 
$\varphi_n(z)$ being a polynomial of order $2n+1$ 
(see \eqref{exp-knun} below).   
$\varphi_n$ has a unique negative root and 
we regard it as a zero of $K_{2n+\frac{3}{2}}$.  
Hence $K_{2n+\frac{3}{2}}$ has $2n+1$ zeros.   
It is known that $K_\nu$ has $2(n+1)$ zeros 
when $2n+\frac{3}{2}<\nu<2n+\frac{7}{2}$.   
It is also well known that the non-real zeros are complex conjugate 
in pairs.   

Recently in \cite{HM}, it has been shown that the zeros of $K_\nu$ 
are obtained as the roots of some algebraic equations 
whose coefficients are explicitly given 
by using $K_\nu$ and $I_\nu$.   
For details, see the equation \eqref{eq-for-zeros} below.   
When $2n+\frac{3}{2}<\nu<2n+\frac{7}{2}$, 
the equations may be taken of order $2(n+1)$.   
Such equations have been already shown in \cite{ZT} 
when $\nu$ is an integer and they coincide in this special case.  

Let $z_\nu$ be a non-real zero of $K_\nu$.   
By the formula 
\begin{equation*}
I_\nu(z)K'_{\nu}(z)-I'_{\nu}(z)K_\nu(z)=-\frac{1}{z}
\end{equation*}
or the uniqueness of the Bessel differential equation, 
we see that $K_\nu'(z_\nu)\ne0$ and 
that $z_\nu$ is (locally) a smooth function in $\nu$ 
by the implicit function theorem.    

The aim of this article is to show the continuity of 
the zeros from the algebraic equations, 
including the continuity at $\nu=2n+\frac{3}{2}$, and 
to present some numerical computations by {\it Mathematica}.   
The following graph shows the behavior of the zeros.  
\begin{figure}[htbp]
\begin{center}
\includegraphics[width=10cm, height=7cm]{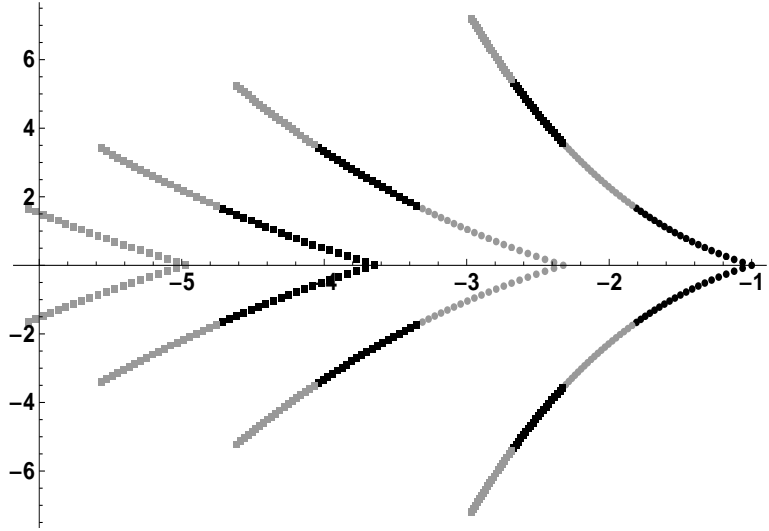}
\end{center}
\caption{Zeros of $K_{\nu}(z)$}
\end{figure}

The unique zero of $K_{\frac{3}{2}}$ is $-1$.   
The two curves from $-1$ described by the black points 
give the two zeros in the case of $\frac{3}{2}<\nu<\frac{7}{2}$.  
The endpoints and the negative value between $-2$ and $-3$ 
found in the graph are the three zeros of $K_{\frac{7}{2}}$.  
The four curves from the zeros of $K_{\frac{7}{2}}$ 
described by the gray points are the zeros in the case of 
$\frac{7}{2}<\nu<\frac{11}{2}$.   
The five zeros of $K_{\frac{11}{2}}$ are seen 
in a similar manner, and so on.   
See also the table in the last part of this article.  
It should be mentioned that we find a similar graph 
for the zeros in \cite{KS}.   

In order to mention the main results, 
we recall some formulae obtained in \cite{HM} 
and prepare some further results.   
For $\nu=\frac{3}{2}$ we have 
\begin{equation*}
K_{\frac{3}{2}}(z)=\sqrt{\frac{\pi}{2z}} \frac{e^{-z}}{z}(1+z).
\end{equation*}
When $\nu>\frac{3}{2}$, we denote the zeros of $K_\nu$ by 
$z_1^{(\nu)},z_2^{(\nu)},..., z_{N(\nu)}^{(\nu)}$.   
Here the number $N(\nu)$ of the zeros of $K_\nu$ 
is equal to $\nu-\frac{1}{2}$ if $\nu-\frac{1}{2}$ is an integer 
and is the even integer closest to $\nu-\frac{1}{2}$ 
otherwise as is mentioned above.  

The basic formula is the following: 
if $\nu-\frac{1}{2}$ is an integer, 
\begin{equation*}
\frac{K_{\nu+1}(w)}{K_\nu(w)}= 1+\frac{2\nu}{w}+
\sum_{j=1}^{N(\nu)}\frac{1}{z_j^{(\nu)}-w}
\end{equation*}
and, if $\nu-\frac{1}{2}$ is not an integer, 
\begin{equation} \label{inv-lap}
\frac{K_{\nu+1}(w)}{K_\nu(w)}= 1+\frac{2\nu}{w}+
\sum_{j=1}^{N(\nu)}\frac{1}{z_j^{(\nu)}-w}+\cos(\pi\nu)
\int_0^\infty \frac{dx}{x(x+w)G_\nu(x)},
\end{equation}
where the function $G_\nu$ is given by 
\begin{align*}
G_\nu(x) & =K_\nu(x)^2+\pi^2I_\nu(x)^2+
2\pi\sin(\pi\nu)K_\nu(x)I_\nu(x) \\
 & = K_\nu(x)^2 + \pi^2 I_\nu(x) I_{-\nu}(x).
\end{align*}
These formulae have been obtained in the course of 
some study on the first hitting times of 
the Bessel diffusion processes.   
By considering the asymptotic expansions of the both hand sides 
of \eqref{inv-lap} as $w\to\infty$, 
we obtain the algebraic equations \eqref{eq-for-zeros} for the zeros 
$z_j^{(\nu)}, j=1,2,...,N(\nu).$   
If we consider the asymptotic behavior as $w\to0$, 
we obtain the equations for the reciprocals, 
which we do not mention in this paper.   

We put 
\begin{equation*}
\nu_n=2n+\frac{3}{2}, \qquad n=0,1,2,...
\end{equation*}
Then, the following is easily seen.   

\medskip

\begin{lemma}
When $\nu=\nu_n,$ $G_\nu$ has a unique positive zero 
and it is the unique solution of $K_{\nu_n}=\pi I_{\nu_n}.$  
When $\nu\ne\nu_n$ for any $n,$ 
$G_\nu$ does not vanish{\rm .}
\end{lemma}

\bigskip

It should be noted that, denoting by $x_n$ 
the unique solution of $K_{\nu_n}=\pi I_{\nu_n}$, 
$-x_n$ is the negative zero of $K_{\nu_n}$, 
which is seen by the formula
\begin{equation*}
K_\nu(e^{m\pi i}z)=e^{-\nu m\pi i}K_\nu(z)-
\pi i \frac{\sin{\nu m\pi}}{\sin(\nu\pi)} I_\nu(z) 
\quad (m\in\mathbf{Z}).
\end{equation*}
Since $K_{\nu_n}$ is explicitly given by 
\begin{equation} \label{exp-knun}
K_{\nu_n}(z)=\sqrt{\frac{\pi}{2z}} e^{-z} \sum_{k=0}^{2n+1}
\frac{(\nu_n,k)}{(2z)^{k}},
\end{equation}
$-x_n$ satisfies 
\begin{equation} \label{eq-for-xn}
\sum_{k=0}^{2n+1} \frac{(\nu_n,2n+1-k)}{2^{2n+1-k}}z^k=0, 
\end{equation}
where $(\nu,0)=1$ and 
\begin{equation*}
(\nu,n)=\frac{(4\nu^2-1^2)(4\nu^2-3^2)\cdots(4\nu^2-(2n-1)^2)}
{n!2^{2n}} = 
\frac{\Gamma(\nu+n+\frac{1}{2})}{n!\Gamma(\nu-n+\frac{1}{2})}.
\end{equation*}
Moreover, by the recurrence relation 
\begin{equation*}
K_\nu(z)-K_{\nu+2}(z)=-\frac{2(\nu+1)}{z}K_{\nu+1}(z),\quad 
I_\nu(z)-I_{\nu+2}(z)=\frac{2(\nu+1)}{z}I_{\nu+1}(z),
\end{equation*}
we can easily show that 
$K_{\nu_n}(x_{n+1})-\pi I_{\nu_n}(x_{n+1})<0$ 
and that the sequence $\{x_n\}_{n=0}^\infty$ is increasing.  

We set 
\begin{equation*}
\alpha_n=-K'_{\nu_n}(x_n)+\pi I'_{\nu_n}(x_n) \qquad \text{and}
\qquad \beta_n=\biggl(-\frac{\pa K_\nu(x_n)}{\pa\nu}+
\pi\frac{\pa I_\nu(x_n)}{\pa\nu}\biggr)\bigg|_{\nu=\nu_n}.
\end{equation*}
Note that $\alpha_n>0$, 
since $K_\nu$ is decreasing and $I_\nu$ is increasing.   
Then, setting $G(x,\nu)=G_\nu(x)$, we have 
\begin{equation*}
G(x_n,\nu_n)=0, \qquad \frac{\pa G}{\pa x}(x_n,\nu_n)=0, \qquad 
\frac{\pa G}{\pa\nu}(x_n,\nu_n)=0
\end{equation*}
and
\begin{equation*}
\frac{\pa^2G}{\pa x^2}(x_n,\nu_n)=2\alpha_n^2, \ 
\frac{\pa^2G}{\pa x\pa\nu}(x_n,\nu_n)=2\alpha_n\beta_n, \ 
\frac{\pa^2G}{\pa\nu^2}(x_n,\nu_n)=2\beta_n^2+
2\pi^2K_{\nu_n}(x_n)^2.
\end{equation*}
Moreover, we can show 
\begin{equation*}
K_{\nu_n}(x_n)=\frac{\pi}{\alpha_n x_n}.
\end{equation*}

Combining the above mentioned formulae, 
we obtain the following.   

\medskip

\begin{lemma} \label{lem-2}
For $m=1,2,...,2n-1,$ it holds that 
\begin{equation*}
\lim_{\nu\to\nu_n\pm0} \cos(\pi\nu) \int_0^\infty
\frac{x^{m-1}}{G(x,\nu)}dx=\pm x_n^{m}.
\end{equation*}
\end{lemma}

\bigskip

The main results are the following.   

\bigskip

\begin{thm} \label{thm-1}
As $\nu\downarrow \nu_n,\ n=0,1,2,...,$ 
two of the zeros of $K_\nu$ converge to $-x_n$ and 
the others to the non-real zeros of $K_{\nu_n}.$ 

\end{thm}

\bigskip

\begin{thm} \label{thm-2}
As $\nu\uparrow\nu_{n},\ n=1,2,...,$ 
each zero of $K_\nu$ converges to a non-real zero 
of $K_{\nu_n}.$
\end{thm}

\bigskip

\noindent{\it Proof of Theorem \ref{thm-1}}\quad 
At first we recall the algebraic equation 
for the zeros of $K_\nu$.   
For this we define $a_k^{(\nu)}, k=0,1,2,...$ inductively by 
\begin{equation} \label{def-ak-nu}
\frac{(\nu+1,m)}{2^m} = \sum_{k=0}^m 
\frac{(\nu,m-k)}{2^{m-k}} a_k^{(\nu)}.
\end{equation}
Moreover we define $\alpha_m^{(\nu)}$ by $\alpha_0^{(\nu)}=1$ and
\begin{equation*}
\alpha_m^{(\nu)}=\frac{1}{m} \sum_{k=1}^m 
\alpha_{m-k}^{(\nu)} \biggl\{ a_{k+1}^{(\nu)} - (-1)^k
\cos(\pi\nu)\int_0^\infty \frac{y^{k-1}}{G_\nu(y)}dy\biggr\}, 
\quad m=1,2,...,2n+1.
\end{equation*}
Then it is shown in \cite{HM} that the zeros of $K_\nu$ are 
the roots of 
\begin{equation} \label{eq-for-zeros} 
\sum_{k=0}^{2n+2} \alpha_{2n+2-k}^{(\nu)} z^k=0
\end{equation}
by computing the asymptotic behavior of the both hands side of 
\eqref{inv-lap} as $w\to\infty$.   

For our purpose we show that, for $m=0,1,...,2n+2$, 
$\lim_{\nu\downarrow\nu_n}\alpha_m^{(\nu)}$ exists and that, 
denoting the limit by $c_m^{(n)}$, 
\begin{equation} \label{eq-factor}
\sum_{k=0}^{2n+2} c_{2n+2-k}^{(n)} z^k = 
(z+x_n) \sum_{k=0}^{2n+1} \frac{(\nu_n,2n+1-k)}{2^{2n+1-k}}z^k.
\end{equation}
Then, since the roots of algebraic equations are continuous 
in the coefficients, we obtain the assertion of the theorem 
from \eqref{eq-for-xn}.   
$-x_n$ is the double root for the polynomial 
in \eqref{eq-factor}. 

For $k=0,1,...,2n+1$, $\nu_n-k+\frac{1}{2}\geqq 1$ and, hence, 
by \eqref{def-ak-nu}, we easily see that $a_k^{(\nu)}$ converges 
as $\nu\downarrow\nu_n$ and that the limit $a_k^{(n)}$ satisfies 
\begin{equation*}
\frac{(\nu_n+1,m)}{2^m} = \sum_{k=0}^m 
\frac{(\nu_n,m-k)}{2^{m-k}}a_k^{(n)},\quad 
m=0,1,...,2n+1.
\end{equation*}
Since $(\nu_n,2n+2)=(\nu_n,2n+3)=(\nu_n+1,2n+3)=0$, 
we see the convergence of $a_{2n+2}^{(\nu)}$ and $a_{2n+3}^{(\nu)}$ 
again by using \eqref{def-ak-nu}.   

Combining the convergence of $a_k^{(\nu)}$ with Lemma \ref{lem-2}, 
we obtain the convergence of $\alpha_m^{(\nu)}$ by induction.   
The limit $c_m^{(n)}$ satisfies 
\begin{equation*}
c_m^{(n)}=\frac{1}{m}\sum_{k=1}^m c_{m-k}^{(n)} 
\{ a_{k+1}^{(n)} - (-1)^k x_n^k \}, \quad 
m=1,2,...,2n+2.
\end{equation*}
From this recurrence relation we get 
\begin{equation} \label{rep-for-cm}
c_m^{(n)}=\frac{(\nu_n,m)}{2^m} + \frac{(\nu_n,m-1)}{2^{m-1}}x_n.
\end{equation}
We can check this by some lengthy computation and 
we omit the details.   

For a proof of the second assertion, we note $c_0^{(n)}=1$ and 
$(\nu_n,2n+2)=0$.   
Then we get from \eqref{rep-for-cm} 
\begin{align*}
\sum_{k=0}^{2n+2}c_{2n+2-k}^{(n)}z^k & = z^{2n+2} + 
\sum_{k=0}^{2n+1} \biggl\{ \frac{(\nu_n,2n+2-k)}{2^{2n+2-k}} + 
\frac{(\nu_n,2n+1-k)}{2^{2n+1-k}}x_n \biggr\} z^k \\
 & = z^{2n+2} + 
\sum_{k=0}^{2n} \frac{(\nu_n,2n+1-k)}{2^{2n+1-k}} z^{k+1} + 
x_n \sum_{k=0}^{2n+1} \frac{(\nu_n,2n+1-k)}{2^{2n+1-k}}z^k \\
 & = \sum_{k=0}^{2n+1} \frac{(\nu_n,2n+1-k)}{2^{2n+1-k}}z^{k+1}+
x_n \sum_{k=0}^{2n+1} \frac{(\nu_n,2n+1-k)}{2^{2n+1-k}}z^k,
\end{align*}
which show the second assertion and Theorem \ref{thm-1}. \qed

\bigskip

\noindent{\it Proof of Theorem \ref{thm-2}}\quad 
We can prove in the same way as in Theorem \ref{thm-1} and 
we only give a sketch.   

When $\nu_{n-1}<\nu<\nu_n$, the zeros of $K_\nu$ are the roots of 
\begin{equation*}
\sum_{k=0}^{2n} \alpha_{2n-k}^{(\nu)}z^k=0.
\end{equation*}
We can show that each $\alpha_m^{(\nu)}, m=0,1,...,2n$ 
converges as $\nu\uparrow\nu_n$ and that, 
denoting the limit by $d_m^{(n)}$, $d_0^{(n)}=1$ and 
\begin{equation*}
d_m^{(n)}=\frac{1}{m} \sum_{k=0}^{m} d_{m-k}^{(n)} 
\{a_{k+1}^{(n)}+(-1)^k x_n^k \}.
\end{equation*}
From this recurrence relation, we obtain 
\begin{equation*}
d_m^{(n)}=\sum_{k=0}^m (-1)^k \frac{(\nu_n,m-k)}{2^{m-k}}x_n^k, 
\quad m=0,1,...,2n
\end{equation*}
and
\begin{equation*}
\sum_{k=0}^{2n+1} \frac{(\nu_n,2n-k+1)}{2^{2n-k+1}}z^k=
(z+x_n)\sum_{k=0}^{2n} d_{2n-k}^{(n)}z^k.
\end{equation*}
The zeros of $K_{\nu_n}$ are the roots of the polynomial 
on the left hand side and 
the limits of the zeros of $K_\nu$ as $\nu\uparrow\nu_n$ 
are the zeros of that on the other side.   
This proves Theorem \ref{thm-2}. \qed

\bigskip

\begin{rem}
It may be worthwhile noting that the algebraic equations 
for the zeros when $\nu_{n-1}<\nu<\nu_{n}$ and 
when $\nu_{n}<\nu<\nu_{n+1}$ are different.   
In fact, letting $z_0$ be one of the zeros of $K_\nu$ 
when $\nu_n<\nu<\nu_{n+1}$, we have 
\begin{equation*}
\sum_{k=0}^{2n+2}\alpha^{(\nu)}_{2n+2-k}z_0^k=0.
\end{equation*}
If $\sum_{k=0}^{2n}\alpha_{2n-k}^{(\nu)}z_0^k=0$, 
we have $\sum_{k=2}^{2n+2}\alpha_{2n+2-k}^{(\nu)}z_0^{k}$
since $z_0\ne0$.   
Then, comparing this equation with the above one, 
we should have that 
$\alpha^{(\nu)}_{2n+2} + \alpha_{2n+1}^{(\nu)} z_0=0$ 
and that $z_0$ is real.   
This is a contradiction.   
\end{rem}

Finally we give a table of the approximate values 
of the zeros of $K_\nu$, 
which are numerically computed by \lq\lq{\it Mathematica}\rq\rq.  
As is mentioned above, 
the algebraic equations for the reciprocals 
of the zeros are also given in \cite{HM}.   
The numerical results for the zeros of the two algebraic equations 
coincide and it gives a good check for our results.   

\newpage 

\def\imag{\mathrm{i}}

\begin{center}
Table of zeros of $K_\nu$
\end{center}
{\tiny 
\begin{center}
\begin{tabular}{|c|c|c|c|c|c|}\hline
 $\nu$ & $\text{zeros}_1$ & $\text{zeros}_2$ & $\text{zeros}_3$ 
& $\text{zeros}_4$ & $\text{zeros}_5$ \\\hline
 $1.5$ & $-1$ & & & \phantom{$-5.03753 \pm 0.0866936 \imag$} &
		     \phantom{$-6.29702$}\\ 
 $1.6$ & $-1.06356 \pm 0.0852232 \imag$ & & & &\\
 $1.7$ & $-1.12292 \pm 0.170806 \imag$ & & & & \\
 $1.8$ & $-1.1787 \pm 0.256725 \imag$ & & & & \\
 $1.9$ & $-1.23139 \pm 0.342957 \imag$ & & & & \\
 $2.0$ & $-1.28137 \pm 0.429485 \imag$ & & & & \\
 $2.1$ & $-1.32896 \pm 0.516291 \imag$ & & & & \\
 $2.2$ & $-1.37442 \pm 0.603361 \imag$ & & & & \\
 $2.3$ & $-1.41795 \pm 0.690682 \imag$ & & & & \\
 $2.4$ & $-1.45976 \pm 0.77824 \imag$ & & & & \\
 $2.5$ & $-1.5 \pm 0.866025 \imag$ & & & & \\
 $2.6$ & $-1.5388 \pm 0.954027 \imag$ & & & & \\
 $2.7$ & $-1.57628 \pm 1.04223 \imag$ & & & & \\
 $2.8$ & $-1.61255 \pm 1.13064 \imag$ & & & & \\
 $2.9$ & $-1.64769 \pm 1.21924 \imag$ & & & & \\
 $3.0$ & $-1.68179 \pm 1.30801 \imag$ & & & & \\
 $3.1$ & $-1.71492 \pm 1.39696 \imag$ & & & & \\
 $3.2$ & $-1.74714 \pm 1.48608 \imag$ & & & & \\
 $3.3$ & $-1.77851 \pm 1.57536 \imag$ & & & & \\
 $3.4$ & $-1.80908 \pm 1.6648 \imag$ & & & & \\\hline
 $3.5$ & $-1.83891 \pm 1.75438 \imag$ & $-2.32219$ & & & \\
 $3.6$ & $-1.86802 \pm 1.84411 \imag$ & $-2.3873 \pm 0.0864217 \imag$ &
	     & & \\
 $3.7$ & $-1.89647 \pm 1.93398 \imag$ & $-2.45036 \pm 0.172895 \imag$ & & & \\
 $3.8$ & $-1.92429 \pm 2.02398 \imag$ & $-2.51152 \pm 0.259427 \imag$ & & & \\
 $3.9$ & $-1.95151 \pm 2.11411 \imag$ & $-2.57092 \pm 0.346026 \imag$ & & & \\
 $4.0$ & $-1.97816 \pm 2.20437 \imag$ & $-2.62867 \pm 0.432697 \imag$ & & & \\
 $4.1$ & $-2.00427 \pm 2.29475 \imag$ & $-2.68489 \pm 0.519443 \imag$ & & & \\
 $4.2$ & $-2.02987 \pm 2.38525 \imag$ & $-2.73967 \pm 0.606267 \imag$ & & & \\
 $4.3$ & $-2.05497 \pm 2.47586 \imag$ & $-2.79309 \pm 0.693173 \imag$ & & & \\
 $4.4$ & $-2.0796 \pm 2.56659 \imag$ & $-2.84525 \pm 0.780161 \imag$ & & & \\
 $4.5$ & $-2.10379 \pm 2.65742 \imag$ & $-2.89621 \pm 0.867234 \imag$ & & & \\
 $4.6$ & $-2.12755 \pm 2.74835 \imag$ & $-2.94604 \pm 0.954392 \imag$ & & & \\
 $4.7$ & $-2.15089 \pm 2.83939 \imag$ & $-2.99479 \pm 1.04164 \imag$ & & & \\
 $4.8$ & $-2.17384 \pm 2.93053 \imag$ & $-3.04252 \pm 1.12897 \imag$ & & & \\
 $4.9$ & $-2.19642 \pm 3.02176 \imag$ & $-3.08929 \pm 1.21638 \imag$ & & & \\
 $5.0$ & $-2.21863 \pm 3.11308 \imag$ & $-3.13513 \pm 1.30388 \imag$ & & & \\
 $5.1$ & $-2.24049 \pm 3.2045 \imag$ & $-3.1801 \pm 1.39147 \imag$ & & & \\
 $5.2$ & $-2.26201 \pm 3.296 \imag$ & $-3.22423 \pm 1.47914 \imag$ & & & \\
 $5.3$ & $-2.28321 \pm 3.38759 \imag$ & $-3.26756 \pm 1.5669 \imag$ & & & \\
 $5.4$ & $-2.30409 \pm 3.47927 \imag$ & $-3.31013 \pm 1.65474 \imag$ & & & \\\hline
 $5.5$ & $-2.32467 \pm 3.57102 \imag$ & $-3.35196 \pm 1.74266 \imag$ &
 $-3.64674$ & & \\
 $5.6$ & $ -2.34497 \pm 3.66286 \imag$ & $ -3.39308 \pm 1.83067 \imag$ &
	     $ -3.71228 \pm 0.0866256 \imag$ & & \\
 $5.7$ & $ -2.36498 \pm 3.75477 \imag$ & $ -3.43354 \pm 1.91875 \imag$ &
	     $ -3.77647 \pm 0.173268 \imag$ & & \\ 
 $5.8$ & $ -2.38471 \pm 3.84676 \imag$ & $ -3.47334 \pm 2.00692 \imag$ &
	     $ -3.83937 \pm 0.259933 \imag$ & & \\
$5.9$ & $ -2.40419 \pm 3.93883 \imag$ & $ -3.51252 \pm 2.09517 \imag$ &
	     $ -3.90105 \pm 0.346624 \imag$ & & \\
$6.0$ & $ -2.4234 \pm 4.03096 \imag$ & $ -3.5511 \pm 2.1835 \imag$ & 
$-3.96156 \pm 0.433345 \imag$ & & \\ 
 $6.1$ & $ -2.44238 \pm 4.12317 \imag$ & $ -3.5891 \pm 2.2719 \imag$ & $ -4.02096 \pm 0.520101 \imag$ & & \\
 $6.2$ & $ -2.46111 \pm 4.21544 \imag$ & $ -3.62654 \pm 2.36038 \imag$ & $ -4.07929 \pm 0.606892 \imag$ & & \\
 $6.3$ & $ -2.47961 \pm 4.30779 \imag$ & $ -3.66344 \pm 2.44894 \imag$ & $ -4.13661 \pm 0.693723 \imag$ & & \\
 $6.4$ & $ -2.49788 \pm 4.4002 \imag$ & $ -3.69983 \pm 2.53757 \imag$ & $ -4.19295 \pm 0.780595 \imag$ & & \\
 $6.5$ & $ -2.51593 \pm 4.49267 \imag$ & $ -3.73571 \pm 2.62627 \imag$ & $ -4.24836 \pm 0.86751 \imag$ & & \\
 $6.6$ & $ -2.53377 \pm 4.58521 \imag$ & $ -3.7711 \pm 2.71505 \imag$ & $ -4.30287 \pm 0.954469 \imag$ & & \\
 $6.7$ & $ -2.55141 \pm 4.67781 \imag$ & $ -3.80603 \pm 2.8039 \imag$ & $ -4.35652 \pm 1.04147 \imag$ & & \\
 $6.8$ & $ -2.56884 \pm 4.77047 \imag$ & $ -3.8405 \pm 2.89282 \imag$ & $ -4.40935 \pm 1.12853 \imag$ & & \\
 $6.9$ & $ -2.58608 \pm 4.86319 \imag$ & $ -3.87453 \pm 2.98181 \imag$ & $ -4.46137 \pm 1.21563 \imag$ & & \\
 $7.0$ & $ -2.60313 \pm 4.95597 \imag$ & $ -3.90813 \pm 3.07087 \imag$ & $ -4.51263 \pm 1.30278 \imag$ & & \\
 $7.1$ & $ -2.61999 \pm 5.0488 \imag$ & $ -3.94131 \pm 3.16 \imag$ & $ -4.56314 \pm 1.38998 \imag$ & & \\
 $7.2$ & $ -2.63667 \pm 5.1417 \imag$ & $ -3.97409 \pm 3.24919 \imag$ & $ -4.61294 \pm 1.47723 \imag$ & & \\
 $7.3$ & $ -2.65318 \pm 5.23464 \imag$ & $ -4.00649 \pm 3.33846 \imag$ & $ -4.66205 \pm 1.56453 \imag$ & & \\
 $7.4$ & $ -2.66951 \pm 5.32764 \imag$ & $ -4.0385 \pm 3.42778 \imag$ &
	     $ -4.71049 \pm 1.65188 \imag$ & & \\\hline
 $7.5$ & $ -2.68568 \pm 5.42069 \imag$ & $ -4.07014 \pm 3.51717\imag$ & 
	     $ -4.75829 \pm 1.73929 \imag$ & $-4.97179$ & \\
$7.6$ & $ -2.70168 \pm 5.5138 \imag$ & $ -4.10142 \pm 3.60663 \imag$ & $ -4.80546 \pm 1.82674 \imag$ & $ -5.03753 \pm 0.0866936
   \imag$ & \\
 $7.7$ & $ -2.71753 \pm 5.60695 \imag$ & $ -4.13236 \pm 3.69615 \imag$ & $ -4.85203 \pm 1.91425 \imag$ & $ -5.10226 \pm 0.173395
   \imag$ & \\
 $7.8$ & $ -2.73322 \pm 5.70016 \imag$ & $ -4.16295 \pm 3.78573 \imag$ & $ -4.89801 \pm 2.00181 \imag$ & $ -5.16602 \pm 0.260108
   \imag$ & \\
 $7.9$ & $ -2.74875 \pm 5.79341 \imag$ & $ -4.19321 \pm 3.87537 \imag$ & $ -4.94342 \pm 2.08942 \imag$ & $ -5.22884 \pm 0.346835
   \imag$ & \\
 $8.0$ & $ -2.76414 \pm 5.88671 \imag$ & $ -4.22315 \pm 3.96507 \imag$ & $ -4.98828 \pm 2.17708 \imag$ & $ -5.29076 \pm 0.433578
   \imag$ & \\
 $8.1$ & $ -2.77939 \pm 5.98006 \imag$ & $ -4.25278 \pm 4.05483 \imag$ & $ -5.0326 \pm 2.2648 \imag$ & $ -5.35181 \pm 0.52034 \imag
  $ & \\
 $8.2$ & $ -2.79449 \pm 6.07346 \imag$ & $ -4.2821 \pm 4.14464 \imag$ & $ -5.07641 \pm 2.35256 \imag$ & $ -5.41202 \pm 0.607123
   \imag$ & \\
 $8.3$ & $ -2.80946 \pm 6.1669 \imag$ & $ -4.31112 \pm 4.23452 \imag$ & $ -5.1197 \pm 2.44038 \imag$ & $ -5.47142 \pm 0.693928 \imag
  $ & \\
 $8.4$ & $ -2.82429 \pm 6.26038 \imag$ & $ -4.33985 \pm 4.32445 \imag$ & $ -5.16251 \pm 2.52825 \imag$ & $ -5.53003 \pm 0.780758
   \imag$ & \\
 $8.5$ & $ -2.83898 \pm 6.35391 \imag$ & $ -4.36829 \pm 4.41444 \imag$ & $ -5.20484 \pm 2.61618 \imag$ & $ -5.58789 \pm 0.867614
   \imag$ & \\
 $8.6$ & $ -2.85355 \pm 6.44748 \imag$ & $ -4.39646 \pm 4.50449 \imag$ & $ -5.2467 \pm 2.70415 \imag$ & $ -5.645 \pm 0.954498 \imag
  $ & \\
 $8.7$ & $ -2.86799 \pm 6.5411 \imag$ & $ -4.42435 \pm 4.59459 \imag$ & $ -5.28812 \pm 2.79217 \imag$ & $ -5.70141 \pm 1.04141 \imag
  $ & \\
 $8.8$ & $ -2.88231 \pm 6.63475 \imag$ & $ -4.45198 \pm 4.68474 \imag$ & $ -5.32909 \pm 2.88025 \imag$ & $ -5.75712 \pm 1.12835
   \imag$ & \\
 $8.9$ & $ -2.89651 \pm 6.72845 \imag$ & $ -4.47935 \pm 4.77495 \imag$ & $ -5.36963 \pm 2.96837 \imag$ & $ -5.81216 \pm 1.21532
   \imag$ & \\
 $9.0$ & $ -2.91058 \pm 6.82219 \imag$ & $ -4.50647 \pm 4.86521 \imag$ & $ -5.40975 \pm 3.05654 \imag$ & $ -5.86655 \pm 1.30233 \imag
  $ & \\
 $9.1$ & $ -2.92454 \pm 6.91597 \imag$ & $ -4.53334 \pm 4.95552 \imag$ & $ -5.44946 \pm 3.14477 \imag$ & $ -5.92031 \pm 1.38936
   \imag$ & \\
 $9.2$ & $ -2.93839 \pm 7.00978 \imag$ & $ -4.55996 \pm 5.04588 \imag$ & $ -5.48877 \pm 3.23304 \imag$ & $ -5.97345 \pm 1.47643
   \imag$ & \\
 $9.3$ & $ -2.95212 \pm 7.10364 \imag$ & $ -4.58635 \pm 5.1363 \imag$ & $ -5.5277 \pm 3.32137 \imag$ & $ -6.026 \pm 1.56354 \imag$ & \\
 $9.4$ & $ -2.96574 \pm 7.19753 \imag$ & $ -4.61251 \pm 5.22676 \imag$ & $ -5.56625 \pm 3.40974 \imag$ & $ -6.07797 \pm 1.65068
   \imag$ & \\\hline
$9.5$ & $ -2.97926 \pm 7.29146 \imag$ & $ -4.63844 \pm 5.31727 \imag$ 
& $ -5.60442 \pm 3.49816 \imag$ & $ -6.12937 \pm 1.73785 \imag$ &
		     $-6.29702$ \\\hline
\end{tabular}
\end{center}
}

\bigskip

\noindent{\bf Acknowledgements}\\
This work is partially supported by Grants-in-Aid for
Scientific Research (B) 22340020, (C) No.24540181 and No.23540183 
of Japan Society for the Promotion of Science (JSPS).

\bigskip


\bigskip

\noindent Yuji Hamana, Department of Mathematics, 
Kumamoto University, 
Kurokami 2-39-1, Kumamoto 860-8555, Japan \\
E-mail: hamana(at)kumamoto-u.ac.jp \\
Hiroyuki Matsumoto, Department of Physics and Mathematics, 
Aoyama Gakuin University, Fuchinobe 5-10-1, Sagamihara 252-5258, 
Japan \\
E-mail: matsu(at)gem.aoyama.ac.jp \\
Tomoyuki Shirai, Institute of Mathematics for Industry, 
Kyushu University, 744, Motooka, Nishi-ku, Fukuoka 819-0395, 
Japan \\
E-mail: shirai(at)imi.kyushu-u.ac.jp

\end{document}